\documentclass{ro_hack}
\usepackage{hyperref}

\DeclareMathOperator*{\argmin}{argmin}
\begin{document}

\title{Another pedagogy for pure-integer Gomory}\thanks{This research was partially supported by NSF grant CMMI--1160915 and ONR grant N00014-14-1-0315}
\runningtitle{Another pedagogy for pure-integer Gomory}
\author{Qi He}\address{IOE Department, University of Michigan, Ann Arbor, MI, USA.\\ \phantom{.}\quad \email{$\{$qihe,jonxlee$\}$@umich.edu}}
\author{Jon Lee}\sameaddress{1}
\date{17 July 2015; revised 16 December 2015}
\begin{abstract}
We present pure-integer Gomory cuts in a way so that they are derived with respect to
 a ``dual form'' pure-integer optimization problem and applied on the standard-form primal side as columns,
 using the primal simplex algorithm.
The input integer problem is not in standard form, and so the cuts are derived a bit differently.
In this manner, we obtain
 a \emph{finitely-terminating} version of pure-integer Gomory cuts that employs the primal
 rather than the dual simplex algorithm.
\end{abstract}
\maketitle
\section*{Introduction}
We assume some familiarity with integer linear optimization;
see \cite{CCZ} for a modern treatment.
We assume  that $A\in \mathbb{Z}_{m\times n}$ has full row-rank $m$, $c\in \mathbb{Z}^n$, and
we consider a pure integer-optimization problem of the ``dual form''
\[
\tag*{($\mbox{D}_\mathcal{I}$)}
\begin{array}{rlcl}
z:= \max & y'b  &      &   \\
      &  y'A  &   \leq  & c'; \\
      &   y  & \in & \mathbb{Z}^m.
\end{array}
\]
The associated continuous relaxation is denoted $\mbox{D}$.

This linear-optimization problem has a non-standard form as a point of departure,
but it is convenient that the dual of the continuous relaxation $\mbox{D}$ has the \emph{standard} ``primal form''
\[
\tag*{($\mbox{P}$)}
\begin{array}{rrcl}
 \min & c'x  &      &   \\
      &  Ax  &   =  & b; \\
      &   x  & \geq & \mathbf{0}.
\end{array}
\]
We note that the only linear programs that we directly solve have the form P,
which is in the appropriate form for treating with the primal simplex algorithm.


For us, the essence of a pure-integer Gomory (or Chv\'atal-Gomory)
cut is to take an inequality $\alpha'x \leq \beta$
that is valid for the continuous relaxation of a pure-integer optimization problem, with $\alpha\in\mathbb{Z}^n$,
and produce the valid cut $\alpha'x \leq \lfloor\beta\rfloor$ for the pure-integer problem.
In the classical presentation of Gomory, the inequality $\alpha'x \leq \beta$ is obtained
by rounding down the left-hand side coefficients of a ``tableau'' equation of a standard-form problem,
which leads to a valid $\alpha'x \leq \beta$ because the variables are non-negative.
In our setting, the integer-constrained variables (on our ``dual side'') are not non-negative,
and our ``tableau'' equations are on the other side (i.e., our ``primal side''),
so we will have to proceed differently.

In \S\ref{sec:classical}, we briefly summarize the classical presentation of pure-integer Gomory cuts.
In \S\ref{sec:neo}, we present a new variation of a pure-integer Gomory cutting-plane
algorithm. Our variation avoids the dual simplex algorithm, and rather precisely fits into the
well-known column-generation framework based on the primal simplex algorithm (commonly used, for example,
for Dantzig-Wolfe decomposition and for the cutting-stock problem; see \cite{LeeLP}, for example).
Furthermore, it turns out there is a certain computational economy
which we will explain.
In \S\ref{sec:example}, we present an example.
In \S\ref{sec:finite}, we present a
finite cutting-plane algorithm based on the cuts of \S\ref{sec:neo}.

An advantage of our set-up and finiteness proof is that it does not rely on the
lexicographical dual simplex method. In the senior author's (30+ years)
experience of of teaching non-doctoral engineering students:
(i) already the lexical primal simplex algorithm is
a topic that many students are challenged to comprehend, but eventually
learn in the context of proving finiteness for the primal simplex algorithm and
for establishing the strong-duality theorem of linear optimization (i.e., the approach in \cite{LeeLP},
for example);
(ii)  a quick explanation of the
dual simplex algorithm is taken as very technical and somewhat mysterious;
(iii) putting these two topics together to explain the lexical dual simplex
algorithm leaves many students behind; so (iv) few students can
then absorb the standard (and quite technical) finiteness proof for classical Gomory cuts,
because it rests on the shaky foundation that they have for the
lexical dual simplex algorithm.

We note that throughout (\S\ref{sec:neo} and \S\ref{sec:finite}),
we deconstruct the derivations and proofs, to completely expose
the movable parts, rather than seeking to make the presentation as short as possible.

\section{Classic Gomory}\label{sec:classical}

In the classical presentation of pure-integer Gomory cutting planes (see \cite{LeeCO}, for example),
we seek to solve a standard-form linear-optimization problem $\mbox{P}$
with the restriction that all variables are integer. Here we assume that $b$ is an integer vector.
From a primal basis $\beta$ for the
standard-form problem (we use $\eta$ for the non-basic indices),
we derive a Gomory cut from any ``source equation''
\[
\tag*{($E_i$)}
x_{\beta_i} + \sum_{j=1}^{n-m} \bar{a}_{i,\eta_j} x_{\eta_j} = \bar{x}_{\beta_i}
\]
having $\bar{x}_{\beta_i}$ non-integer. The cut is simply
\[
x_{\beta_i} + \sum_{j=1}^{n-m} \left\lfloor\bar{a}_{i,\eta_j}\right\rfloor x_{\eta_j} \leq \left\lfloor\bar{x}_{\beta_i}\right\rfloor,
\]
which is clearly violated by $x=\bar{x}$. Introducing a non-negative slack variable $x_k$,
we get an equation
\[
x_{\beta_i} + \sum_{j=1}^{n-m} \left\lfloor\bar{a}_{i,\eta_j}\right\rfloor x_{\eta_j} + x_k = \left\lfloor\bar{x}_{\beta_i}\right\rfloor,
\]
which, subtracting the equation $E_i$, can be introduced at the current stage as
\[
 \sum_{j=1}^{n-m} \left(\left\lfloor\bar{a}_{i,\eta_j}\right\rfloor -  \bar{a}_{i,\eta_j}\right)
  x_{\eta_j} + x_k = \left\lfloor\bar{x}_{\beta_i}\right\rfloor - \bar{x}_{\beta_i}.
\]
This new variable $x_k$ is an additional basic variable, but it has a negative value
$x_k = \left\lfloor\bar{x}_{\beta_i}\right\rfloor - \bar{x}_{\beta_i}$.
Naturally, we proceed to re-optimize by the dual simplex algorithm, seeking to regain primal feasibility while maintaining
dual feasibility. As established by Gomory (see \cite{Gomory}; also see \cite{LeeCO} for a presentation more closely following
the notation used here), this can be realized as a \emph{finite} algorithm
by:
\begin{itemize}
\item[(i)] introducing an integer objective variable $x_0$ and
associated equation $x_0-c'x=0$; here we note that it is important that
$c$ is an integer vector, and so $x_0$ is an integer on the set of feasible integer solutions of
$\mbox{P}$; moreover, the additional equation implies an additional basic variable,
which we take as $x_0$ and deem it to be the first basic variable: that is, $\beta \leftarrow (0,\beta)$,
now an ordered list of $m+1$ basic indices from $\{0,1,2,\ldots,n\}$,
\item[(ii)] always choosing a source equation $E_i$ with least $i$ among those with
$\bar{x}_{\beta_i}$ non-integer; here we stress the importance of
the objective variable $x_0$ having index $0$ and being the first basic variable at the outset,
\item[(iii)]
sequentially numbering added slack variables $x_{n+1},x_{n+2},\ldots$,
\item[(iv)] re-solving each linear-optimization problem after a cut via the \emph{lexicographic (i.e., epsilon-perturbed) dual simplex algorithm} (see \cite{Lemke}).
\end{itemize}

Gomory did say in \cite{Gomory}:
\begin{quote}
\emph{``In these proofs we will use the lexicographical dual simplex method described in Section 7. It
is not implied that this simplex method be used in practice or that it is necessary to the proof.
It is simply that its use in the proof has reduced the original rather long and tedious proofs
to relatively simple ones.''}
\end{quote}

This proof has endured in all presentations
that we know of (e.g,
see
\cite[pp. 215--6]{CCZ};
\cite[pp. 165--7]{LeeCO};
\cite[pp. 372--3]{NW};
\cite[pp. 285--7]{PR};
\cite[pp. 121--3]{SM};
\cite[pp. 354--8]{Schrijver}),
and we do not know Gomory's ``original rather long and tedious proofs''.
Incidentally, many of these published proofs are lacking a bit in complete rigor,
including the one of the the second author of the present paper (see \cite{LeeCO}).
A clear unfortunate aspect of the proof is its delicate set up.

The classical way of doing Gomory uses the dual simplex algorithm
 because the cut-generation methodology seems wedded to a standard form for the integer problem that
 we wish to solve.
In what follows, we derive cuts a bit differently so that the dual of the
continuous relaxation of the integer problem that
 we wish to solve is in standard form.
In this way, we simply add columns to a standard-form problem and naturally re-optimize
via the primal simplex algorithm. Though just a bit more complicated in its derivation than the classical Gomory approach, our method can be presented and implemented in a unified and simple manner with other column-generation algorithms based on the primal simplex algorithm  (in particular,
Dantzig-Wolfe decomposition, the cutting-stock algorithm, and even
a presentation of \emph{branch-and-bound}; see \cite{LeeLP}). Moreover, the actual calculations are quite straightforward to carry out (see the example in \S\ref{sec:example}). Finally, we wish to point out and emphasize that
in carrying out the primal simplex algorithm for P, every basis has precisely $m$ elements, even as we add columns. This is in sharp contrast to the
classical Gomory approach, where each cut increments the number of basis elements (as well as appends a slack variable).
If classical Gomory were to be applied to our
formulation $\mbox{D}_\mathcal{I}$ (which has $m$ variables
in $n$ inequalities), putting it into standard form
would give us a problem with $2m+n$ (non-negative) variables in
$n$ equations. So (dual simplex algorithm) bases would have size $n$
\emph{and} would grow as we add cuts. Because of this, the matrix algebra of each pivot-step
in our approach is simpler.

We note that \cite{LeeWiegele2015} addresses extending our approach to the mixed-integer case. Another direction that could be explored is how to lift inequalities to strengthen them (see \cite{DeyRichard}, for example).

\section{Gomory another way}\label{sec:neo}

Let us return to approaching the pure integer-optimization problem $\mbox{D}_\mathcal{I}$.
In what follows, we refer to $\mbox{D}$ (the continuous relaxation of $\mbox{D}_\mathcal{I}$) as the \emph{dual}
and $\mbox{P}$ (the dual of  $\mbox{D}$) as the \emph{primal}.
Let $\beta$ be any basis for $\mbox{P}$. The associated dual basic solution (for the continuous relaxation
$\mbox{D}$) is $\bar{y}':=c_\beta'A_\beta^{-1}$. Suppose that $\bar{y}_i$ is not an integer. Our goal is to derive a valid cut for $\mbox{D}_\mathcal{I}$
that is violated by $\bar{y}$.

Let
\[
\tilde{b}:=\mathbf{e}^i +  A_\beta r,
\]
where $r\in\mathbb{Z}^m$, and
$\mathbf{e}^i$ denotes the $i$-th standard unit
vector in $\mathbb{R}^m$.
Note that by construction, $\tilde{b}\in\mathbb{Z}^m$.

\begin{thrm} \label{lem:nonint}
$\bar{y}'\tilde{b}$
is not an integer, and so $y'\tilde{b} \leq \lfloor \bar{y}' \tilde{b} \rfloor$
cuts off $\bar{y}$.
\end{thrm}
\begin{proof}
$\bar{y}'\tilde{b}= \bar{y}'
(\mathbf{e}^i + A_\beta r)=
\bar{y}_i + (c_\beta' A_\beta^{-1}) A_\beta r
= \underbrace{\bar{y}_i}_{\notin\mathbb{Z}} + \underbrace{c_\beta'r}_{\in\mathbb{Z}}$.
\end{proof}



At this point, we have an inequality
$y'\tilde{b} \leq \lfloor \bar{y}' \tilde{b} \rfloor$
which cuts off $\bar{y}$, but we have not established its
validity for $\mbox{D}_\mathcal{I}$.

Let $H_{\cdot i}:=A_\beta^{-1}\mathbf{e}^i$,  the $i$-th column of $A_\beta^{-1}$. Now let
\[
w:= H_{\cdot i} + r.
\]
Clearly we can choose
$r\in\mathbb{Z}^m$
so that $w\geq \mathbf{0}$;
we simply choose
$r\in\mathbb{Z}^m$ so that
\begin{equation}\label{rcondition}
r_k \geq - \lfloor h_{ki} \rfloor, \mbox{ for } k=1,\ldots,m.
\end{equation}

\begin{thrm}
Choosing $r\in\mathbb{Z}^m$ satisfying
(\ref{rcondition}),
we have that
$y'\tilde{b} \leq \lfloor \bar{y}' \tilde{b} \rfloor$
is valid for $\mbox{D}_\mathcal{I}$.
\end{thrm}
\begin{proof}
Because $w\geq 0$ and $y'A\leq c'$, we have the validity of
\[
y'A_\beta (A_\beta^{-1}\mathbf{e}^i + r) \leq c_\beta'(A_\beta^{-1}\mathbf{e}^i + r),
\]
even for the continuous relaxation D of $\mbox{D}_\mathcal{I}$.
Simplifying this, we have
\[
y'(\mathbf{e}^i + A_\beta r) \leq \bar{y}_i + c_\beta'r.
\]
The left-hand side is clearly $y'\tilde{b}$, and the right-hand side
is
\[
\bar{y}_i + c_\beta'r=\bar{y}_i + \bar{y}'A_\beta r =
\bar{y}'(\mathbf{e}^i + A_\beta r) =
\bar{y}' \tilde{b}.
\]
So we have that
$y'\tilde{b} \leq \bar{y}' \tilde{b}$
is valid even for $\mbox{D}$.
%
Finally, observing that $\tilde{b}\in\mathbb{Z}^m$ and $y$ is constrained to
be in $\mathbb{Z}^m$ for $\mbox{D}_\mathcal{I}$,
we can round down the right-hand side and get the result.
\end{proof}

So, given any  non-integer basic dual solution $\bar{y}$,
we have a way to produce a valid inequality for $\mbox{D}_\mathcal{I}$ that cuts it off.
This cut for $\mbox{D}_\mathcal{I}$ is used as a column for $\mbox{P}$:
the column is $\tilde{b}$ with objective coefficient $\lfloor \bar{y}' \tilde{b} \rfloor$.
Taking $\beta$ to be an optimal basis for $\mbox{P}$,
the new variable corresponding to this column is the unique variable eligible to enter the basis
in the context of the primal simplex algorithm applied to $\mbox{P}$ --- the reduced cost is precisely
\[
\bar{y}' \tilde{b} - \lfloor \bar{y}' \tilde{b} \rfloor <0.
\]

\begin{bsrvtn}
The new column for $A$ is $\tilde{b}$ which is integer.
The new objective coefficient for $c$ is $\lfloor \bar{y}' \tilde{b} \rfloor$ which is an integer.
So the original assumption that $A$ and $c$ are integer is maintained, and we can repeat.
In this way, we get a legitimate cutting-plane framework for $\mbox{D}_\mathcal{I}$
 --- though we emphasize that we do our computations as column generation
with respect to $\mbox{P}$.
\end{bsrvtn}

There is clearly a lot of flexibility in
how $r$ can be chosen.
Next, we demonstrate that in a very concrete sense, it is always best to choose
a minimal $r\in\mathbb{Z}^m$ satisfying (\ref{rcondition}).

\begin{thrm}
Let
$r\in\mathbb{Z}^m$ be defined by
\begin{equation}\label{rcondition_minimal}
r_k = - \lfloor h_{ki} \rfloor, \mbox{ for } k=1,\ldots,m,
\end{equation}
and suppose  that $\hat{r}\in\mathbb{Z}^m$ satisfies
$r\leq \hat{r}$. Then the cut determined
by $r$ dominates the cut determined by $\hat{r}$.
\end{thrm}

\begin{proof}
It is easy to check that our cut can be re-expressed as
\[
y_i \leq \lfloor \bar{y}_i\rfloor + \left( c_\beta' -y'A_\beta\right)r.
\]
Noting that $c_\beta' -y'A_\beta\geq \mathbf{0}$ for all $y$ that are feasible
for $\mbox{D}$, we see that the strongest inequality
is obtained by choosing $r\in\mathbb{Z}^m$ to be minimal.
\end{proof}

\section{Example}\label{sec:example}

In this section, we present an example which illustrates the
simplicity of the calculations. Throughout, we choose $r\in\mathbb{Z}^m$
to be minimal, as defined in \ref{rcondition_minimal}.

Let
\[
A=\left(
  \begin{array}{ccccc}
    7 & 8 & -1 & 1 & 3 \\
    5 & 6 & -1 & 2 & 1 \\
  \end{array}
\right), \quad
b=\left(
    \begin{array}{c}
      26 \\
      19 \\
    \end{array}
  \right)
\]
\[
\mbox{and } c' = \left(
      \begin{array}{ccccc}
        126 & 141 & -10 & 5 & 67 \\
      \end{array}
    \right).
\]
So, the integer program $\mbox{D}_\mathcal{I}$ which we seek to solve is defined by five inequalities in the two variables $y_1$ and $y_2$.
For the basis of $\mbox{P}$, $\beta=(1,2)$, we have
\[
A_\beta=\left(\begin{array}{cc}
          7 & 8 \\
          5 & 6
        \end{array}
        \right),
\mbox{ and hence }
A_\beta^{-1}
=
\left(
  \begin{array}{cc}
    3 & -4 \\
    -5/2 & 7/2 \\
  \end{array}
\right).
\]
It is easy to check that for this choice of basis, we have
\[
\bar{x}_\beta=\left(
                \begin{array}{c}
                  2 \\
                  3/2 \\
                \end{array}
              \right),
\]
and for the non-basis $\eta=\{3,4,5,6\}$, we have
$\bar{c}_\eta'=\left(\begin{array}{ccc}
                5 & 1/2 & 1
                              \end{array}\right)$,
which are both non-negative, and so this basis is optimal for
$\mbox{P}$. The associated dual basic solution is
\[
\bar{y}'=\left(
           \begin{array}{cc}
             51/2 & -21/2 \\
           \end{array}
         \right), \mbox{ and the objective value is } z=463~1/2.
\]

Because both $\bar{y}_1$ and $\bar{y}_2$ are not integer, we
can derive a cut for $\mbox{D}_\mathcal{I}$ from either.
Recalling the procedure, for any fraction $\bar{y}_i$,
we start with the $i$-th column $H_{\cdot i}$ of $H:=A_\beta^{-1}$,
and we get a new $A_{\cdot j} := \mathbf{e}^i +  A_\beta r$. That is,
\[
H_{\cdot 1}= \left(\!\!\!
               \begin{array}{c}
                 3 \\
                 -5/2 \\
               \end{array}\!\!\!
             \right)
             \Rightarrow
             r=\left(\!\!\!
                 \begin{array}{c}
                   -3 \\
                   3 \\
                 \end{array}\!\!\!
               \right)
                \Rightarrow
                \tilde{b}=\left(
                            \begin{array}{c}
                              1 \\
                              0 \\
                            \end{array}
                          \right)
                          +
\left(\begin{array}{cc}
          7 & 8 \\
          5 & 6
        \end{array}
        \right)\left(
                 \begin{array}{c}
                   -3 \\
                   3 \\
                 \end{array}
               \right)
               = \left(
                   \begin{array}{c}
                     4 \\
                     3 \\
                   \end{array}
                 \right) =: A_{\cdot 6}
\]

\[
H_{\cdot 2}= \left(\!\!\!
               \begin{array}{c}
                 -4 \\
                 7/2 \\
               \end{array}\!\!\!
             \right)
                          \Rightarrow
             r=\left(\!\!\!
                 \begin{array}{c}
                   4 \\
                   -3 \\
                 \end{array}\!\!\!
               \right)
                \Rightarrow
                \tilde{b}=\left(
                            \begin{array}{c}
                              0 \\
                              1 \\
                            \end{array}
                          \right)
                          +
\left(\begin{array}{cc}
          7 & 8 \\
          5 & 6
        \end{array}
        \right)\left(
                 \begin{array}{c}
                   4 \\
                   -3 \\
                 \end{array}
               \right)
               = \left(
                   \begin{array}{c}
                     4 \\
                     3 \\
                   \end{array}
                 \right).
\]
In fact, for this iteration of this example, we get the same cut for either
choice of $i$. To calculate the right-hand side of the cut, we have
\[
\bar{y}'\tilde{b} = \left(
           \begin{array}{cc}
             51/2 & -21/2 \\
           \end{array}
         \right) \left(
                   \begin{array}{c}
                     4 \\
                     3 \\
                   \end{array}
                 \right) = 70~1/2,
\]
so the cut for $\mbox{D}_\mathcal{I}$ is
\[
4y_1 + 3y_2 \leq 70.
\]

Now, we do our simplex-method calculations with respect to $\mbox{P}$. The new column for $\mbox{P}$ is
$A_{\cdot 6}$ (above) with objective coefficient $c_6:= 70$.

Following the ratio test of the primal simplex algorithm,
when index 6 enters the basis,
index 2 leaves the basis, and so the new basis is $\beta=(1,6)$, with
\[
A_\beta=\left(\begin{array}{cc}
          7 & 4 \\
          5 & 3
        \end{array}
        \right),
\]
with objective value 462, a decrease.
At this point, index 5 has a negative reduced cost, and index 1 leaves the basis. So we now have $\beta=(5,6)$,
which turns out to be optimal for the current $\mbox{P}$.
We have
\[
\bar{y}'=\left(
           \begin{array}{cc}
             131/5 & -58/5 \\
           \end{array}
         \right), \mbox{ and the objective value is }  z=460~4/5.
\]
We observe that the objective function has decreased, but unfortunately
both $\bar{y}_1$ and $\bar{y}_1$ are not integers. So we must continue.
We have
\[
A_\beta=\left(\begin{array}{cc}
          3 & 4 \\
          1 & 3
        \end{array}
        \right),
\mbox{ and hence }
A_\beta^{-1}
=
\left(
  \begin{array}{cc}
    3/5 & -4/5 \\
    -1/5 & 3/5 \\
  \end{array}
\right).
\]

We observe that the objective function has decreased,
but because both $\bar{y}_1$ and $\bar{y}_2$ are not integers, we
can again derive a cut for $\mbox{D}_\mathcal{I}$ from either.
We calculate
\[
H_{\cdot 1}= \left(\!\!\!
               \begin{array}{c}
                 3/5 \\
                 -1/5 \\
               \end{array}\!\!\!
             \right)
             \Rightarrow
             r=\left(
                 \begin{array}{c}
                   0 \\
                   1 \\
                 \end{array}
               \right)
                \Rightarrow
                \tilde{b}=\left(
                            \begin{array}{c}
                              1 \\
                              0 \\
                            \end{array}
                          \right)
                          +
\left(\begin{array}{cc}
          3 & 4 \\
          1 & 3
        \end{array}
        \right)\left(
                 \begin{array}{c}
                   0 \\
                   1 \\
                 \end{array}
               \right)
               = \left(
                   \begin{array}{c}
                     5 \\
                     3 \\
                   \end{array}
                 \right)=:A_{\cdot 7}
\]

\[
H_{\cdot 2}= \left(\!\!\!
               \begin{array}{c}
                 -4/5 \\
                 3/5 \\
               \end{array}\!\!\!
             \right)
                          \Rightarrow
             r=\left(
                 \begin{array}{c}
                   1 \\
                   0 \\
                 \end{array}
               \right)
                \Rightarrow
                \tilde{b}=\left(
                            \begin{array}{c}
                              0 \\
                              1 \\
                            \end{array}
                          \right)
                          +
\left(\begin{array}{cc}
          3 & 4 \\
          1 & 3
        \end{array}
        \right)\left(
                 \begin{array}{c}
                   1 \\
                   0 \\
                 \end{array}
               \right)
               = \left(
                   \begin{array}{c}
                     3 \\
                     2 \\
                   \end{array}
                 \right)=:A_{\cdot 8}~ .
\]
Correspondingly, we have $\bar{y}'A_{\cdot 7}= 96~1/5$ and $\bar{y}'A_{\cdot 8}= 55~2/5$, giving us $c_7:= 96$ and $c_8 := 55$. So, we have two possible cuts for
$\mbox{D}_\mathcal{I}$:
\[
5y_1 + 3y_2 \leq 96 \mbox{ and } 3y_1 + 2y_2 \leq 55.
\]

Choosing to incorporate both as columns for $\mbox{P}$, and
letting index 8 enter the basis, index 5 leaves (according to the primal-simplex ratio test),
and it turns out that we reach an optimal basis $\beta=(8,6)$
after this single pivot. At this point, we have
\[
\bar{y}'=\left(
           \begin{array}{cc}
             25 & -10 \\
           \end{array}
         \right), \mbox{ and the objective value is }  z=460.
\]
Not only has the objective decreased, but now all of the $\bar{y}_i$
are integers, so we have an optimal solution for $\mbox{D}_\mathcal{I}$.

We wish to emphasize that to take this example with $n=5$ inequalities in $m=2$ unrestricted variables and put it into standard form, we would end up
with $2m+n=9$ variables and $n=5$ equations. So, the initial basis for
applying the classical Gomory algorithm would have $n=5$ elements, and
subsequent bases after cuts would be even larger. In contrast, our bases have
$m=2$ elements throughout, thus making the matrix algebra less burdensome.

\section{Finite convergence}\label{sec:finite}

To make a finitely-converging algorithm,
we amend our set-up a bit:
\begin{itemize}
\item[(i)] we assume that the objective vector $b$ is integer, and we
move the objective function to the constraints;
\item[(ii)] \emph{after this}, we lexicographically perturb the resulting objective
function.
\end{itemize}
So, we arrive at
\[
\tag*{($\mbox{D}^{\epsilon}_\mathcal{I}$)}
\begin{array}{rlclcl}
 \max & y_{0} & +    & y'\vec{\epsilon}_{\scriptscriptstyle[1,m]}    &      &   \\
      & y_{0} & -    & y'b                 & \leq & 0;\\
      &       &      & y'A                 & \leq & c'; \\
      & \multicolumn{3}{l}{y_{0}  \in  \mathbb{Z};}  &      &     \\
      &       &      & \multicolumn{3}{l}{y  \in  \mathbb{Z}^m,}
\end{array}
\]
where  $\vec{\epsilon}_{\scriptscriptstyle[i,j]}:=(\epsilon^i,\epsilon^{i+1},...,\epsilon^{j})'$, and $\epsilon$ is
treated as an arbitrarily small positive \emph{indeterminate} --- we wish to emphasize that
we do not give $\epsilon$ a real value, rather we incorporate it symbolically.
We note that if $(y_0,y')$ is optimal for
$\mbox{D}^{\epsilon}_\mathcal{I}$, then $y$ is a lexically-maximum solution of
$\mbox{D}_\mathcal{I}$; that is, $y$ is optimal for $\mbox{D}_\mathcal{I}$,
and it is lexically maximum (among all optimal solutions)
under the total ordering of basic dual solutions induced by $\sum_{i=1}^m \epsilon^i y_i$.

The dual of the continuous relaxation of $\mbox{D}^{\epsilon}_\mathcal{I}$ is the rhs-perturbed
primal problem
\[
\tag*{($\mbox{P}^{\epsilon}$)}
\begin{array}{rrcrcl}
 \min &       &      & c'x    &      &   \\
      & x_0   &      &        & =    & 1;\\
      &-bx_0  & +    & Ax     & =    & \vec{\epsilon}_{\scriptscriptstyle[1,m]}; \\
      & \multicolumn{3}{l}{x_0  \geq  0;}  &      &     \\
      &       &      & \multicolumn{3}{l}{x  \geq  \mathbf{0},}
\end{array}
\]

Next, we observe that $\mbox{D}^{\epsilon}_\mathcal{I}$
is a special case of
\[
\tag*{($\mbox{lex-D}_\mathcal{I}$)}
\begin{array}{rlcl}
z:= \max & y'\vec{\epsilon}_{\scriptscriptstyle[0,m-1]}   &      &   \\
      &  y'A  &   \leq  & c'; \\
      &   y  & \in & \mathbb{Z}^m,
\end{array}
\]
which has as the dual of its continuous relaxation the rhs-perturbed primal problem
\[
\tag*{($\mbox{lex-P}$)}
\begin{array}{rrcl}
 \min & c'x  &      &   \\
      &  Ax  &   =  &  \vec{\epsilon}_{\scriptscriptstyle[0,m-1]}; \\
      &   x  & \geq & \mathbf{0}.
\end{array}
\]
So, in what follows, we focus on $\mbox{lex-D}_\mathcal{I}$ and $\mbox{lex-P}$.

\subsection{First pivot after a new column}

The primal simplex algorithm applied to the
\emph{non-degenerate} $\mbox{lex-P}$
produces a sequence of dual solutions $\bar{y}'$ with \emph{decreasing}
objective value $\bar{y}'\vec{\epsilon}_{\scriptscriptstyle[0,m-1]}$.
This can be interpreted as a \emph{lexically decreasing} sequence of $\bar{y}$.
We wish to emphasize that after we add a new column to $\mbox{lex-P}$,
on the next pivot (and of course subsequent ones), the
basic dual solution $\bar{y}$ lexically decreases. We want to show more.

\begin{lmm}
If we derive a column from an $i$ for which $\bar{y}_i$ is fractional
(in the manner of \S\ref{sec:neo}), append this column to $\mbox{lex-P}$,
and then make a single primal-simplex pivot, say with the $l$-th basic variable leaving the basis,
then after the pivot the new dual solution is
\[
\bar{\bar{y}}=  \bar{y} + \frac{\lfloor\bar{y}_i\rfloor - \bar{y}_i}{h_{li}+r_l}H_{l\cdot},
\]
where $H_{l\cdot}$ is the $l$-th row of $A_{\beta}^{-1}$.
\end{lmm}

\begin{proof}
This is basic simplex-algorithm stuff. $\bar{\bar{y}}$ is
just $\bar{y}$ plus a multiple $\Delta$ of the $l$-th for of $A_\beta^{-1}$.
The  reduced cost of the entering variable, which starts
at $\lfloor\bar{y}_i\rfloor - \bar{y}_i$
 will become zero (because it becomes basic) after the
pivot. So
\[
\left(\lfloor\bar{y}_i\rfloor - \bar{y}_i\right) - \Delta \left( h_{li}+r_l \right) =0,
\]
which implies that
\[
\Delta = \frac{\lfloor\bar{y}_i\rfloor - \bar{y}_i}{ h_{li}+r_l}.
\]
\end{proof}

\begin{crllr}\label{cor:minimal_step}
If we derive a column from an $i$ for which $\bar{y}_i$ is fractional
(in the manner of \S\ref{sec:neo}), choosing $r\in\mathbb{Z}^m$
to be minimal (i.e., satisfies (\ref{rcondition_minimal})),
append this column to $\mbox{lex-P}$,
and then make a single primal-simplex pivot,
then after the pivot, either $(\bar{\bar{y}}_1,\ldots,\bar{\bar{y}}_{i-1})$ is a lexical decrease
relative to $(\bar{y}_1,\ldots,\bar{y}_{i-1})$
or $\bar{\bar{y}}_i\leq \lfloor\bar{y}_i\rfloor$.
\end{crllr}

\begin{proof}
A primal pivot implies that we observe the usual ratio test to maintain primal feasibility. This
amounts to choosing
\[
l :=  \argmin_{l ~:~ h_{li} + r_l > 0} \left\{
\frac{H_{l\cdot} \vec{\epsilon}_{\scriptscriptstyle[0,m-1]}}{h_{li} + r_l}
\right\}.
\]
Also, we have
\[
\bar{\bar{y}}_i=  \bar{y}_i + \frac{\overbrace{\lfloor\bar{y}_i\rfloor - \bar{y}_i}^{<0}}{\underbrace{h_{li} + r_l}_{>0}}h_{li}.
\]
Assume that  $(\bar{\bar{y}}_1,\ldots,\bar{\bar{y}}_{i-1})$ is not a lexical decrease
relative to $(\bar{y}_1,\ldots,\bar{y}_{i-1})$.
Because $\bar{\bar{y}}$ is lexically less than $\bar{y}$, we then must have  $h_{li}\geq 0$.

\[
\bar{\bar{y}}_i=  \bar{y}_i + \frac{\lfloor\bar{y}_i\rfloor - \bar{y}_i}{h_{li}+r_l}h_{li}
= \bar{y}_i +
   \underbrace{\left(\lfloor\bar{y}_i\rfloor - \bar{y}_i\right)}_{\displaystyle <0}
   \underbrace{\left(\frac{h_{li}}{h_{li}+r_l}\right)}_{
  \begin{array}{c}
    \geq 1?
  \end{array}
   }
\leq \bar{y}_i + \left(\lfloor\bar{y}_i\rfloor - \bar{y}_i\right)
= \lfloor\bar{y}_i\rfloor.
\]
To finish the proof, we need to justify
\[
\tag*{($\Phi$)} \frac{h_{li}}{h_{li}+r_l} \geq 1.
\]
A sufficient condition for $\Phi$  to hold is $r_l\leq 0$ and $h_{li}>0$.
Taking
 $r$ to be minimal, we have $h_{li} + r_l=h_{li}-\lfloor h_{li}\rfloor>0$
 which, together with $h_{li}\geq 0$, implies that $h_{li}>0$  and $r_l=-\lfloor h_{li}\rfloor \leq 0$
\end{proof}

\begin{bsrvtn}
We note that we are using the fact that we choose $r$ to be minimal
to get $\Phi$ to hold.
 However, it is not necessary that
we choose $r\in\mathbb{Z}^m$ to be minimal for the conclusion of Corollary \ref{cor:minimal_step}
to hold.
We simply need to have $r_l\leq 0$ and $h_{li}>0$ to ensure that
$\Phi$ holds.
\end{bsrvtn}

\subsection{A finite column-generation algorithm for pure integer-linear optimization}

Next, we specify a finitely-converging algorithm for $\mbox{lex-D}_\mathcal{I}$.
We assume that the feasible region of the continuous relaxation $\mbox{D}$ of
$\mbox{D}_\mathcal{I}$  is non-empty and bounded.
Because of how we reformulate $\mbox{D}_\mathcal{I}$
as $\mbox{lex-D}_\mathcal{I}$, we have that the feasible region of
the associated continuous relaxation $\mbox{lex-D}$
is non-empty and bounded.

\vfill\eject

\begin{center}
\underline{\bf Algorithm 1: Column-generation for pure integer-linear optimization}
\end{center}
\begin{enumerate}
\item[(0)] Assume that the feasible region of $\mbox{lex-D}$
is non-empty and bounded.
Start with the basic feasible optimal solution of $\mbox{lex-P}$ (obtained in any manner).
\item Let $\bar{y}$ be the associated dual basic solution. 
If $\bar{y}_i \in\mathbb{Z}$ for all $i\in\mathcal{I}$, then STOP: $\bar{y}$ solves $\mbox{lex-D}_\mathcal{I}$.\label{alg:STOP}
\item Otherwise, choose the \emph{minimum} $i\in\mathcal{I}$ for which $\bar{y}_i \notin\mathbb{Z}$. Related to this $i$,
construct a new variable (and associated column and objective coefficient) for $\mbox{lex-P}$
in the manner of \S\ref{sec:neo}, choosing $r$ to be \emph{minimal}. Solve this new version of $\mbox{lex-P}$,
starting from the current (primal feasible) basis, employing the primal simplex algorithm.
\begin{enumerate}
\item If this new version of $\mbox{lex-P}$ is unbounded, then STOP: $\mbox{lex-D}_\mathcal{I}$ is infeasible.
\item Otherwise, GOTO step \ref{alg:STOP}. \label{alg:cut}
\end{enumerate}
\end{enumerate}

%
%
%
%
%
%

\begin{thrm}
Algorithm 1 terminates in a finite number of iterations with either an optimal solution of $\mbox{lex-D}_\mathcal{I}$
or a proof that $\mbox{lex-D}_\mathcal{I}$ is infeasible.
\end{thrm}

\begin{proof}
It is clear from well-known facts about linear optimization that if the algorithm stops,
then the conclusions asserted by the algorithm are correct. So our task is to
demonstrate that the algorithm terminates in a finite number of iterations.

Consider the full sequence of dual solutions $\bar{y}^t$ ($t=1,2,\ldots$) visited during the algorithm.
We refer to every dual solution after every \emph{pivot} (of the primal-simplex algorithm), over all visits to step
\ref{alg:cut}.
This sequence is lexically decreasing at every (primal-simplex) pivot.
We claim that after a finite number of iterations of Algorithm 1,
$\bar{y}^t$ is an integer vector upon reaching step \ref{alg:STOP}, whereupon the algorithm stops.
If not, let
$j$ be the least index for which $\bar{y}_j$ does not become and remain constant (and integer) after a finite number of pivots

 Choose an iteration $T$ where $\bar{y}^T$ of step \ref{alg:STOP} has
$\bar{y}^T_k$  constant (and integer) for all $k<j$ and all
 subsequent pivots.
 Consider the infinite (non-increasing) sequence $\mathcal{S}_1:=\bar{y}^T_j, \bar{y}^{T+1}_j, \bar{y}^{T+2}_j,\cdots$.
 By the choice of $j$, this sequence has an infinite strictly decreasing subsequence $\mathcal{S}_2$.
 By the boundedness assumption, this subsequence has an infinite strictly decreasing subsequence $\mathcal{S}_3$
 of fractional values that are between some pair of successive integers.
 By Corollary \ref{cor:minimal_step}, between any two visits to step \ref{alg:STOP}
 with $\bar{y}_j$ fractional, there is at least one integer between these fractional values.
 Therefore,
 $\mathcal{S}_3$  corresponds to pivots in the same visit to step \ref{alg:cut}.
 But this contradicts the fact that the lexicographic primal simplex algorithm converges in a finite number of iterations.
\end{proof}

\bibliographystyle{amsplain}
\bibliography{gp}

\providecommand{\bysame}{\leavevmode\hbox to3em{\hrulefill}\thinspace}
\providecommand{\MR}{\relax\ifhmode\unskip\space\fi MR }
\providecommand{\MRhref}[2]{%
  \href{http://www.ams.org/mathscinet-getitem?mr=#1}{#2}
}
\providecommand{\href}[2]{#2}
\begin{thebibliography}{10}

\bibitem{CCZ}
Michele Conforti, G{\'e}rard Cornu{\'e}jols, and Giacomo Zambelli,
  \emph{Integer programming}, Graduate Texts in Mathematics, vol. 271,
  Springer, 2014.

\bibitem{DeyRichard}
Santanu~S. Dey and Jean-Philippe Richard, \emph{Linear-programming-based
  lifting and its application to primal cutting-plane algorithms}, INFORMS
  Journal on Computing \textbf{21} (2009), no.~1, 137--150.

\bibitem{Gomory}
Ralph~E. Gomory, \emph{An algorithm for integer solutions to linear programs},
  Recent advances in mathematical programming, McGraw-Hill, New York, 1963,
  pp.~269--302.

\bibitem{LeeCO}
Jon Lee, \emph{A first course in combinatorial optimization}, Cambridge Texts
  in Applied Mathematics, Cambridge University Press, Cambridge, 2004.

\bibitem{LeeLP}
Jon\vspace{0mm} Lee, \emph{A first course in linear optimization ({S}econd
  edition, version 2.1)}, Reex Press, 2013--5,
  \url{https://github.com/jon77lee/JLee_LinearOptimizationBook}.

\bibitem{LeeWiegele2015}
Jon Lee\vspace{0mm} and Angelika Wiegele, \emph{Another pedagogy for
  mixed-integer gomory}, Tech. report, 2015.

\bibitem{Lemke}
Carlton~E. Lemke, \emph{The dual method of solving the linear programming
  problem}, Naval Research Logistics Quarterly \textbf{1} (1954), 36--47.

\bibitem{NW}
George~L. Nemhauser and Laurence~A. Wolsey, \emph{Integer and combinatorial
  optimization}, Wiley-Interscience Series in Discrete Mathematics and
  Optimization, John Wiley \& Sons, Inc., New York, 1988, A Wiley-Interscience
  Publication.

\bibitem{PR}
R.~Gary Parker and Ronald~L. Rardin, \emph{Discrete optimization}, Computer
  Science and Scientific Computing, Academic Press, Inc., Boston, MA, 1988.

\bibitem{SM}
Harvey~M. Salkin and Kamlesh Mathur, \emph{Foundations of integer programming},
  North-Holland Publishing Co., New York, 1989.

\bibitem{Schrijver}
Alexander Schrijver, \emph{Theory of linear and integer programming},
  Wiley-Interscience Series in Discrete Mathematics, John Wiley \& Sons, Ltd.,
  Chichester, 1986, A Wiley-Interscience Publication.

\end{thebibliography}
\end{document}